\documentclass[10pt,leqno]{amsart}
\usepackage{graphicx}
\usepackage{indentfirst,csquotes}

\usepackage{amssymb,amsthm,amsmath}
\usepackage{xcolor,paralist,hyperref,titlesec,fancyhdr,etoolbox}
\newtheorem{theorem}{Theorem}[]
\newtheorem{definition}[theorem]{Definition}

\newtheorem{lemma}[theorem]{Lemma}

\titleformat{\section}[display]{\normalfont\huge\bfseries\centering}{\centering\chaptertitlename\thechapter}{10pt}{\Large}
\titlespacing*{\section}{0pt}{0ex}{0ex}

\hypersetup{ colorlinks=true, linkcolor=black, filecolor=black, urlcolor=black }
\def\proof{\noindent {\it Proof. $\, $}}

\def\?#1{#1\nobreak\discretionary{}{\hbox{$\mathsurround=0pt #1$}}{}}

\begin{document}
\title{Planarity ranks of modular varieties of semigroups} %%%%%%%%%%%%
\author[Solomatin D.V.]{Solomatin Denis Vladimirovich}
\date{\today}
\address{Omsk State Pedagogical University, Tukhachevsky embankment, 14}
\email{solomatin\_dv@omgpu.ru}
\maketitle

\let\thefootnote\relax

\begin{abstract}
By the planarity rank of a semigroup variety we mean the largest number of generators of a free semigroup of a variety with respect to which the semigroup admits a planar Cayley graph. Since the time when L.M.~Martynov formulated the problem of describing the planarity ranks of semigroup varieties, many specific results have been obtained in this direction. A modular variety of semigroups is a variety of semigroups with a modular lattice of subvarieties. In this paper, we calculate the exact values of the planarity ranks of an infinite countable set of all possible modular varieties of semigroups. It turns out that these values do not exceed 3. Machine calculations are mostly used in the proof. Prover9 and Mace4 are used to check the equalities of elements of free semigroups of varieties defined by a large number of identities. To prove the non-planarity of graphs, the Pontryagin--Kuratovsky criterion is used, and the Colin de Verdi\`ere invariant is indirectly used to justify planarity. 
\end{abstract} %%%%%%%%%

\bigskip

\noindent 
%*********************************************************
\section*{Introduction}%The introduction and conclusion are not numbered
%*********************************************************

The study of the possibility of representing semigroups using graphs actively occupies the minds of contemporaries~\cite{0}, but we consider such a representation in the context of calculating the planarity ranks of semigroup varieties.

Recall~\cite{1} that a modular variety of semigroups is a variety of semigroups with a modular lattice of subvarieties, in other words, with a lattice in which each pair of elements $a\ ,\ b\ \in \ L$ is modular, that is, the modularity law holds --- the quasi-identity: $\forall \ x\ \in \ L\ :\ x\ \leqslant \ b\ \?\Rightarrow \ x\ \vee \ (\ a\ \wedge \ b\ )\ =\ (\ x\ \vee \ a\ )\ \wedge \ b$. Any non-modular lattice contains the five-element pentagon $N_5$ as a sublattice.

The study of modular varieties of semigroups has several important aspects. From a theoretical point of view, modular varieties of semigroups help to understand the fundamental properties of algebraic structures. They allow us to study various identities and structural properties of semigroups, which is important for the general theory of semigroups, as it helps to construct various generalizations and classifications. Studying modular varieties allows us to classify semigroups by their structural properties and identities. This helps to create a more general and systematic picture of various algebraic objects. Modular varieties of semigroups also have applied significance. The results of research in this area can be applied in other areas of mathematics and theoretical computer science. For example, they can be useful in the theory of automata and formal languages. Modular varieties of semigroups play an important role in the theory of automata and formal languages due to the following aspects. Firstly, modular varieties of semigroups help to classify formal languages by their structural properties. This allows us to better understand which languages can be recognized by certain types of automata. Secondly, using modular varieties allows us to optimize finite automata, which makes them more efficient for recognizing and processing languages. Third, modular varieties of semigroups are used to analyze regular languages and their properties. This is important for developing algorithms that can work effectively with these languages. Fourth, they are also used in automata synthesis, which allows one to create automata that can recognize more complex languages or perform more complex tasks.

Let us give a general example of how the modular variety of semigroups can be used in automata theory and formal languages to optimize finite automata. Suppose that there is a finite automaton that recognizes a certain language. This automaton can be represented as a semigroup, where the states of the automaton and the transitions between them form the structure of a semigroup. Consider an automaton $A$ with a set of states $Q$ and a set of input symbols $\mathrm{\Sigma }$. The transitions of the automaton can be described by the function $\delta :Q\times \mathrm{\Sigma }\to Q$. The semigroup $S$ consists of all possible transitions of the automaton. Using the modular variety of semigroups, one can classify the transitions and states of the automaton. This allows one to identify equivalent states and transitions that can be combined to optimize the automaton. Using the theory of modular varieties, one can minimize the automaton by reducing the number of states and transitions. This is done by combining equivalent states and transitions, which leads to a more compact and efficient automaton. Suppose there is an automaton that recognizes the language of all strings containing an even number of symbols $a$. Then the states of the automaton can be divided into two classes: states where the number of $a$ is even, and states where the number of $a$ is odd. Using the modular variety of semigroups, one can combine all states with an even number of $a$ into one state, and all states with an odd number of $a$ into another state. This will lead to the minimization of the automaton to two states, which significantly simplifies its structure and improves performance. Thus, modular varieties of semigroups help optimize finite automata, making them more efficient for recognizing and processing formal languages. Modular varieties of semigroups have long been characterized in various ways. In the following lemma we present one of the known characterizations in the language of identities.

\begin{lemma}[\cite{1}, Theorem 1 (equational version) as a result of \cite{2}] A variety of semigroups is modular if and only if it satisfies one of the following systems of identities (where $n$ is a natural number):

($m1_n$) : $xy\ =\ {\left(xy\right)}^{n+1}\ $; 

($m2_n$) : $xy\ =x^{n+1}y$, ${\left(xy\right)}^{n+1}=xy^{n+1}$, $xyzt=xyx^nzt$; 

($m3_n$) : $xy=xy^{n+1}$, ${\left(xy\right)}^{n+1}=x^{n+1}y$, $xyzt=xyt^nzt$; 

($m4_{n,\pi }$) : $x_1x_2x_3x_4=x_{1\pi }x_{2\pi }x_{3\pi }x_{4\pi }$, $x^2y=xyx=yx^2=x^{n+2}y$ ($\pi \in {\prod}_1$); 

($m5_{\pi }$) : $x_1x_2x_3x_4=x_{1\pi }x_{2\pi }x_{3\pi }x_{4\pi }$, $x^2y=xyx=yx^2$, $x^3yz=xy^3z$, $x^6=x^7$ ($\pi \in {\prod}_1$); 

($m6_{\pi }$) : $x_1x_2x_3x_4=x_{1\pi }x_{2\pi }x_{3\pi }x_{4\pi }$, $x^2y=xyx=yx^2$, $x^2y^2z=xy^2z^2$ ($\pi \in {\prod}_1$); 

($m7_{\pi }$) : $x_1x_2x_3x_4=x_{1\pi }x_{2\pi }x_{3\pi }x_{4\pi }$, $x^2y=xyx=yx^2$, $x^3yz=xy^2z^2$ ($\pi \in {\prod}_1$); 

($m8_{\pi }$) : $x_1x_2x_3x_4=x_{1\pi }x_{2\pi }x_{3\pi }x_{4\pi }$, $x^2y=xyx$, $xy^2=yx^2$ ($\pi \in {\prod}_1$); 

($m9_{\pi }$) : $x_1x_2x_3x_4=x_{1\pi }x_{2\pi }x_{3\pi }x_{4\pi }$, $x^2y=y^2x$, $xy^2=yxy$ ($\pi \in {\prod}_1$); 

($m{10}_{\pi }$) : $x_1x_2x_3x_4=x_{1\pi }x_{2\pi }x_{3\pi }x_{4\pi }$, $x^2y=yxy$, $xy^2=yx^2$ ($\pi \in {\prod}_1$); 

($m{11}_{\pi }$) : $x_1x_2x_3x_4=x_{1\pi }x_{2\pi }x_{3\pi }x_{4\pi }$, $x^2y=y^2x$, $xy^2=xyx$ ($\pi \in {\prod}_1$); 

($m{12}_{\pi }$) : $x_1x_2x_3x_4=x_{1\pi }x_{2\pi }x_{3\pi }x_{4\pi }$, $x^2y=xy^2$, $xyx=yxy$ ($\pi \in {\prod}_1$); 

($m{13}_{\pi }$) : $x_1x_2x_3x_4=x_{1\pi }x_{2\pi }x_{3\pi }x_{4\pi }$, $x^2y=yxy=yx^2$ ($\pi \in {\prod}_1$); 
 
($m{14}_{\pi }$) : $x_1x_2x_3x_4=x_{1\pi }x_{2\pi }x_{3\pi }x_{4\pi }$, $x^2y=xyx=xy^2$, $x^4y=yx^4$ ($\pi \in {\prod}_1$); 

($m{15}_{\pi }$) : $x_1x_2x_3x_4=x_{1\pi }x_{2\pi }x_{3\pi }x_{4\pi }$, $x^2y=yxy=xy^2$, $x^4y=yx^4$ ($\pi \in {\prod}_1$); 

($m{16}_{\pi }$) : $x_1x_2x_3x_4=x_{1\pi }x_{2\pi }x_{3\pi }x_{4\pi }$, $x^2y=y^2x$, $xyx=x^2yx$, $x^3y=yx^3$ ($\pi \in {\prod}_1$); 

($m{17}_{\pi }$) : $x_1x_2x_3x_4=x_{1\pi }x_{2\pi }x_{3\pi }x_{4\pi }$, $xy^2=yx^2$, $xyx=xyx^2$, $x^3y=yx^3$ ($\pi \in {\prod}_1$); 

($m{18}_{\pi }$) : $x_1x_2x_3x_4=x_{1\pi }x_{2\pi }x_{3\pi }x_{4\pi }$, $x^2y=x^3y$, $xyx=yxy$, $x^3y=yx^3$ ($\pi \in {\prod}_1$); 

($m{19}_{\pi }$) : $x_1x_2x_3x_4=x_{1\pi }x_{2\pi }x_{3\pi }x_{4\pi }$, $xy^2=xy^3$, $xyx=yxy$, $x^3y=yx^3$ ($\pi \in {\prod}_1$); 

($m{20}_{\pi }$) : $x_1x_2x_3x_4=x_{1\pi }x_{2\pi }x_{3\pi }x_{4\pi }$, $x^2y=yx^2$, $xyx=yxy$ ($\pi \in {\prod}_2$); 

($m21$) : $xyzt=tzyx$, $x^2y=yx^2$, $xyx=yxy$, $x^2yz=y^2zx$;

($m22$) : $xyzt=tzyx$, $x^2y=yx^2$, $xyx=yxy$, $x^2yz=yzyx$; 

($m23$) : $xyzt=tzyx$, $x^2y=yx^2$, $xyx=yxy$, $xyxz=yzyx$; 

($m{24}_{\pi }$) : $x_1x_2x_3x_4=x_{1\pi }x_{2\pi }x_{3\pi }x_{4\pi }$, $x^2y=x^3y$, $xy^2=yx^2$, $x^3y=yx^3$ ($\pi \in {\prod}_3$); 

($m{25}_{\pi }$) : $x_1x_2x_3x_4=x_{1\pi }x_{2\pi }x_{3\pi }x_{4\pi }$, $x^2y=y^2x$, $xy^2=xy^3$, $x^3y=yx^3$ ($\pi \in {\prod}_3$); 

($m26$) : $xyzt=ztxy$, $x^2y=x^3y$, $xy^2=yx^2$, $xyxz=yxzx$; 

($m27$) : $xyzt=ztxy$, $x^2y=x^3y$, $xy^2=yx^2$, $xyxz=yxyz$; 

($m28$) : $xyzt=ztxy$, $x^2y=x^3y$, $xy^2=yx^2$, $xyzy=xzyz$; 

($m29$) : $xyzt=ztxy$, $x^2y=y^2x$, $xy^2=xy^3$, $xyxz=yxzx$; 

($m30$) : $xyzt=ztxy$, $x^2y=y^2x$, $xy^2=xy^3$, $xyxz=yxyz$; 

($m31$) : $xyzt=ztxy$, $x^2y=y^2x$, $xy^2=xy^3$, $xyzy=xzyz$; 

($m32$) : $xyzt=tzyx$, $x^2y=x^3y$, $xy^2=yx^2$, $xyxz=xyzx$; 

($m33$) : $xyzt=tzyx$, $x^2y=x^3y$, $xy^2=yx^2$, $xyxz=yxzy$; 

($m34$) : $xyzt=tzyx$, $x^2y=x^3y$, $xy^2=yx^2$, $xyxz=zxyz$; 

($m35$) : $xyzt=tzyx$, $x^2y=x^3y$, $xy^2=yx^2$, $xyxz=yzyx$; 

($m36$) : $xyzt=tzyx$, $x^2y=x^3y$, $xy^2=yx^2$, $xyzx=yxzy$; 

($m37$) : $xyzt=tzyx$, $x^2y=y^2x$, $xy^2=xy^3$, $xyxz=xyzx$; 

($m38$) : $xyzt=tzyx$, $x^2y=y^2x$, $xy^2=xy^3$, $xyxz=yxzy$; 

($m39$) : $xyzt=tzyx$, $x^2y=y^2x$, $xy^2=xy^3$, $xyxz=zxyz$; 

($m40$) : $xyzt=tzyx$, $x^2y=y^2x$, $xy^2=xy^3$, $xyxz=yzyx$; 

($m41$) : $xyzt=tzyx$, $x^2y=y^2x$, $xy^2=xy^3$, $xyzx=yxzy$; 

($m42$) : $xyzt=yxtz$, $x^2y=y^2x$, $xy^2={\left(xy\right)}^2$; 

($m43$) : $xyzt=yxtz$, $x^2y={\left(xy\right)}^2$, $xy^2=yx^2$; 

($m44$) : $xyzt=ztxy$, $x^2y=y^2x$, $xy^2={\left(xy\right)}^2$, $xyxz=yxzx$; 

($m45$) : $xyzt=ztxy$, $x^2y={\left(xy\right)}^2$, $xy^2=yx^2$, $xyxz=yxzx$; 

($m46$) : $xyzt=tzyx$, $x^2y=y^2x$, $xy^2={\left(xy\right)}^2$, $xyzx=yxzy$; 

($m47$) : $xyzt=tzyx$, $x^2y={\left(xy\right)}^2$, $xy^2=yx^2$, $xyzx=yxzy$. 

where 
\[{\prod}_1=\{(123),\ (124),\ (134),\ (234),\ (12)(34),\ (13)(24),\ (14)(23)\};\] 
\[{\prod}_2=\{(123),\ (124),\ (134),\ (234),\ (12)(34),\ (13)(24)\};\] 
\[{\prod}_3=\{(123),\ (124),\ (134),\ (234),\ (12)(34)\}.\] 
\end{lemma} % Lemma

Closely related to the Cayley graphs of semigroups, planar automata, in particular cellular automata, also have many practical applications and are actively studied for several reasons. For example, they are used to model various natural and artificial systems, such as the spread of fire, population growth, the spread of diseases, and even crowd behavior. Planar automata help to study fundamental questions in the theory of computation, such as the universality and complexity of algorithms. Some cellular automata are used in cryptographic applications to create complex and break-resistant ciphers, and can also be used to generate pseudo-random numbers, which is important for various applications, including statistics and modeling.

$\,$
$\,$
%******************************************************
\section*{Main Result}%1
%******************************************************

\begin{definition} Let $S$ be a semigroup and $X$ be the set of elements forming it. Then, the \textit{Cayley graph} $Cay(S,X)$ of a semigroup $S$ with respect to the set of its forming elements $X$ is a directed multigraph with labeled edges whose vertex set coincides with $S$, and whose vertex $a$ is connected to vertex $b$ by an arc starting at vertex $a\in S$, ending at vertex $b\in S$ and labeled by element $x \in X$ if and only if the equality $ax=b$ holds in the semigroup $S$.
\end{definition}

\begin{definition} \textit{The simplified Cayley graph} $SCay(S,X)$ of $Cay(S,X)$ is an simple graph obtained from the original graph $Cay(S,X)$ by removing loops, labels, and replacing all arcs connecting the same pairs of vertices with a single edge connecting the same vertices.
\end{definition}

\begin{definition} According to the well-known Pontryagin--Kuratowski theorem, an ordinary graph is \textit{planar} if and only if it does not contain subgraphs homeomorphic to the complete graph $K_5$ of order five or to the complete bipartite graph $K_{3,3}$ containing three vertices in each of the parts.
\end{definition}

\begin{definition} \textit{A semigroup} $S$ \textit{admits a planar Cayley graph} $Cay(S,X)$ if, with respect to some minimal set of indecomposable generators of $X$, the simplified Cayley graph $SCay(S,X)$ is a planar graph.
\end{definition}

\begin{definition} Let $\textbf{V}=\text{var}\{ \sigma \} $ be a variety of semigroups defined by the system of identities $\sigma$, and $F_n(\textbf{V})$ be a free semigroup in the variety $\textbf{V}$ generated by $n$ generators. Then, the \textit{planarity rank of the manifold} $\textbf{V}$ is a natural number $r_{\pi}(\textbf{V})$ such that for $n\leq r_{\pi}(\textbf{V})$ the semigroup $F_n(\textbf{V})$ admits a planar Cayley graph, and for $n > r_{\pi}(\textbf{V})$ the semigroup $F_n(\textbf{V})$ does not admit a planar Cayley graph.
\end{definition}

The following theorem solves L.M.~Martynov~\cite{mart}'s problem of describing planarity ranks in the class of modular varieties of semigroups and contains a complete list of planarity ranks of these varieties.

\begin{theorem}
$r_{\pi }\left(\mathrm{var}\{m1_1\}\right)=2$, $r_{\pi }\left(\mathrm{var}\{m1_n\}\right)=1$, when $n\ge 2$.

$r_{\pi }\left(\mathrm{var}\{m2_1\}\right)=r_{\pi }\left(\mathrm{var}\{m2_2\}\right)=2$, $r_{\pi }\left(\mathrm{var}\{m2_n\}\right)=1$, when $n\ge 3$.

$r_{\pi }\left(\mathrm{var}\{m3_1\}\right)=2$, $r_{\pi }\left(\mathrm{var}\{m3_n\}\right)=1$, when $n\ge 2$.

$r_{\pi }\left(\mathrm{var}\{m4_{1\pi }\}\right)=2$, $r_{\pi }\left(\mathrm{var}\{m4_{n\pi }\}\right)=1$, when $n\ge 2$, $\pi \in {\prod}_1$.

$r_{\pi }\left(\mathrm{var}\{m5_{\pi }\}\right)=r_{\pi }\left(\mathrm{var}\{m6_{\pi }\}\right)=r_{\pi }\left(\mathrm{var}\{m7_{\pi }\}\right)=1$, when $\pi \in {\prod}_1$.

$r_{\pi }\left(\mathrm{var}\{m8_{\pi }\}\right)=2$, when $\pi \in {\prod}_1$.

$r_{\pi }\left(\mathrm{var}\{m9_{\pi }\}\right)=r_{\pi }\left(\mathrm{var}\{m{10}_{\pi }\}\right)=1$, when $\pi \in {\prod}_1$.

$r_{\pi }\left(\mathrm{var}\{m{11}_{\pi }\}\right)=r_{\pi }\left(\mathrm{var}\{m{12}_{\pi }\}\right)=r_{\pi }\left(\mathrm{var}\{m{13}_{\pi }\}\right)=2$, when $\pi \in {\prod}_1$.

$r_{\pi }\left(\mathrm{var}\{m{14}_{\pi }\}\right)=3$, when $\pi \in {\prod}_1$.

$r_{\pi }\left(\mathrm{var}\{mN_{\pi }\}\right)=2$, when $\pi \in {\prod}_1$, where $N=15\div 19$.

$r_{\pi }\left(\mathrm{var}\{m{20}_{\pi }\}\right)=r_{\pi }\left(\mathrm{var}\{m21\}\right)=r_{\pi }\left(\mathrm{var}\{m22\}\right)=r_{\pi }\left(\mathrm{var}\{m23\}\right)=1$, when $\pi \in {\prod}_2$.

$r_{\pi }\left(\mathrm{var}\{m{24}_{\pi }\}\right)=r_{\pi }\left(\mathrm{var}\{m{25}_{\pi }\}\right)=2$, when $\pi \in {\prod}_3$.

$r_{\pi }\left(\mathrm{var}\{mN\}\right)=2$, where $N=26\div 42$.

$r_{\pi }\left(\mathrm{var}\{m43\}\right)=r_{\pi }\left(\mathrm{var}\{m45\}\right)=r_{\pi }\left(\mathrm{var}\{m47\}\right)=1$.

$r_{\pi }\left(\mathrm{var}\{m44\}\right)=r_{\pi }\left(\mathrm{var}\{m46\}\right)=2$.
\end{theorem} % Theorem

\proof The main problem that has to be solved when constructing Cayley graphs in the general case is the algorithmically unsolvable problem of the equality of words in a semigroup. For example, a request for a detailed proof of the equality $abbaba=ababba$ of elements of a free semigroup of a variety defined by an identity ($m1_2)$ in Prover9 has the form:

\noindent formulas(assumptions).(x*y)*z=x*(y*z).x1*x2=((((x1*x2)*x1)*x2)*x1)*x2.end\_of\_list.

\noindent formulas(goals).((((c1*c2)*c2)*c1)*c2)*c1=((((c1*c2)*c1)*c2)*c2)*c1.end\_of\_list.

The resulting detailed proof takes 128 transformations. The proof of this equality is briefly as follows:

1. $xyyxyx = xyxyyx$; equality to be proven

2. $(xy)z = x(yz)$; original

3. $xy = xyxyxy$; original

4. $xyxyxy = xy$; applied 3

5. $xyz = xyxyxyz$; applied 2 to 4

6. $xyxyxyz = xyz$; applied 5

7. $xyzu = xyzxyzxyzu$; applied 2 to 6

8. $xyz xyz xyz u = xyz u$; applied 7

9. $xyxxyxxy = xyxyxy$; applied 4 to 8

10. $xy x xy x xy = xy$; applied 4 to 9

11. $x (yz xyz xyz u) x x (yz xyz xyz u) x x (yz u = x (yz xyz xyz u)$; applied 8 to 10

12. $xyz xyz xyz u x xyz xyz xyz u x xyz u = xyz xyz xyz u$; applied 2 to 11

13. $xyz xyz xyz u x xyz u x xyz u = xyz xyz xyz u$; applied 8 to 12

14. $xyzxyzxyzuxxyzuxxyzu = xyzxyzxyzu$; applied 2 to 13

15. $xyz u x xyz u x xyz u = xyz xyz xyz u$; applied 8 to 14

16. $xyz u x xyz u x xyz u = xyz u$; applied 8 to 15

17. $xyxyx = xyyxyxyyxyx$; applied 16 to 6

18. $(xy)yxyxyyxyx = x y x yx$; applied 17 by grouping $(x yxyx y) yxyx y yxyx = x y x y x$

19. $xyxxyxyxy = xyxyxxyxy$; applied 18 to 8

20. $xyxxy = xyxyxxyxy$; applied 4 to 19

21. $xyxyxxyxy = xyxxy$; applied 20

22. $xyxyyx = xyyxyx$; applied 21 to 6.

We will not dwell on each of the equality and inequality that arises along the way, we will only note that the reader will be able to check their truth independently using such special-purpose software tools as Prover9 and Mace4, respectively. To do this, three parallel threads are launched: Prover9 in an attempt to prove equality, at the same time it is launched in an attempt to prove inequality, and the third thread launches Mace4 in an attempt to find a counterexample of equality, to prove inequality. The decision is made depending on which thread stops faster with the found answer. Launching a separate Mace4 thread in an attempt to find a counterexample of inequality does not make sense, since a trivial counterexample for inequality, in the form of a specific equality $0=0$, is usually not difficult to find, but its existence does not prove equality in the general case.

Based on Lemma 1, we consider each of the varieties listed in it. If the condition of the theorem being proved is satisfied, then a planar embedding of the Cayley graph of the free semigroup of the manifold under consideration with the corresponding number of generators will be given. As the number of generators increases, forbidden configurations of the type $K_5$ or $K_{3,3}$ are found in the simplified Cayley graph of this semigroup, therefore, according to the Pontryagin--Kuratovsky theorem, the graph will not be planar. Thus, based on the definition, we calculate all planarity ranks. In this case, a short notation of $x-y$ route fragments means that in the set of generators there is an element $z$ such that $xz=y$ or $yz=x$.

$m1$:

A planar embedding of the graph $SCay(F_2(\mathrm{var}\{m1_1\}),\{a,b\})$ is shown in Fig. 1.

\begin{center}
\includegraphics{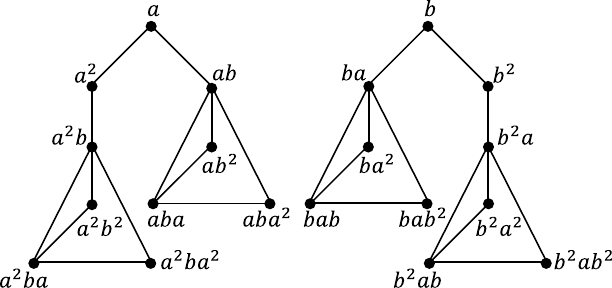}\\
{\footnotesize Fig. 1. Planar embedding of the graph $SCay(F_2(\mathrm{var}\{m1_1\}),\{a,b\})$.}
\end{center}

The presence of a subgraph homeomorphic to the graph $K_{3,3}$ in the simplified Cayley graph $SCay(F_3(\mathrm{var}\{m1_1\}),\{a,b,c\})$ is due to the existence of the following pairwise disjoint routes between the vertices of the sets $\left\{abcb^2,abca,abc^2\right\}$ and $\left\{abcab,abc,abcb\right\}$: $abcb^2-abcb^2a-abcbac-abcaba-abcab$; $abcb^2-abc$; $abcb^2-abcb$; $abca-abcab$; $abca-abc$; $abca-abcac-abcacb-abcba-abcb$; $abc^2-abc^2a-abcab$; $abc^2-abc$; $abc^2-abcb$.

Since in each semigroup variety defined by the identity $x=x^n$ the equality $xy={(xy)}^n$ holds, then by~\cite[Theorem]{3} the planarity rank of the variety $\mathrm{var}\{m1_m\}$ does not exceed 1 for $m\ge 2$. In particular, the graph $Cay(F_3(\mathrm{var}\{m1_2\}),\{x_1,x_2\})$ contains the subgraph shown in the figure~\cite[Fig.2]{3}, the simplified graph of which is homeomorphic to the graph $K_{3,3}$.

$m2$:

The planar embedding of the graph $SCay(F_2(\mathrm{var}\{m2_1\}),\{a,b\})$ and planar embedding of the scheme of one of the two isomorphic connected components of the $SCay(F_2(\mathrm{var}\{m2_2\}),\{a,b\})$ are shown in Fig.2.

In the graph $SCay(F_3(\mathrm{var}\{m2_1\}),\{a,b,c\})$, a subgraph homeomorphic to the graph $K_5$ is found on the following pairwise disjoint routes between vertices that are elements of the set $\{abc, abac, abca, abcb, abc^2\}$: $abc-ab-aba-abac$; $abc-abca$; $abc-abcb$; $abc-abc^2$; $abac-abca$; $abac-abcb$; $abac-abc^2$; $abca-abcb$; $abca-abc^2$; $abcb-abc^2$. Therefore, the Cayley graph of the free semigroup of the manifold $\mathrm{var}\{m2_1\}$ with respect to three or more generators is not planar. In $SCay(F_3(\mathrm{var}\{m2_2\}),\{a,b,c\})$, a subgraph homeomorphic to the graph $K_{3,3}$ can be similarly found on routes between vertices from the set $\{ab,aba^2,ab^3\}$ in one part and $\{aba,ab^2,abc\}$ in the other: $ab-aba$; $ab-ab^2$; $ab-abc$; $aba^2-aba$; $aba^2-ab^2$; $aba^2$$-aba^2c-abca-abc$; $ab^3-aba$; $ab^3-ab^2$; $ab^3-abc$. And then in $SCay(F_2(\mathrm{var}\{m2_n\}),\{a,b\})$, for $n\ge 3$, a subgraph homeomorphic to the graph $K_{3,3}$ is found, in particular for $n=3$ on the routes between the vertices from the set $\{ab^2,aba^2,aba\}$ in one part and $\{aba^3,ab^3,{\left(ab\right)}^2b\}$ in the other: $ab^2-aba^3$; $ab^2-ab^3$; $ab^2-ab^2a-ab^2a^2-ab^2a^2b-ab^2a^2ba-aba^2{\left(ba\right)}^2-{\left(ab\right)}^2ba^2-{\left(ab\right)}^2ba-{\left(ab\right)}^2b$; $aba^2-aba^3$; $aba^2-aba^2b-aba^2ba-a{\left(ba^2\right)}^2-a{\left(ba^2\right)}^2a-ab^3ab^2-ab^3ab-ab^3a-ab^3$; $aba^2-{\left(ab\right)}^2b^2-{\left(ab\right)}^2b$; $aba-aba^3$; $aba-ab^4-ab^3$; $aba-{\left(ab\right)}^2-{\left(ab\right)}^2b$. Intuitively, the value $r_{\pi }(\mathrm{var}\{m2_n\})=1$, for $n\ge 3$, is due to the fact that in each manifold $\mathrm{var}\{m2_n\}$, for $n\ge 3$, the identity $x^2=x^{n+2}$ holds, for $n\ge 3$ leading to the appearance of a simple cycle of three or more elements, as when the identity $x=x^n$ holds, for $n\ge 4$, for the manifold generated by which the planarity rank is equal to 1 [3, Theorem].

\begin{center}
\includegraphics{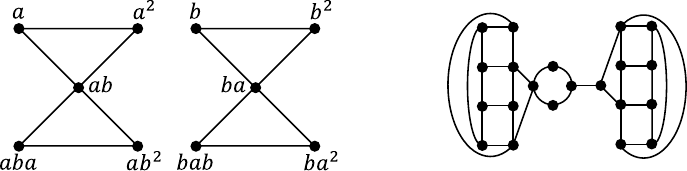}\\
{\footnotesize Fig. 2. Planar embedding of a pair of isomorphic components forming the graph $SCay(F_2(\mathrm{var}\{m2_1\}),\{a,b\})$ [left] and a diagram of the planar embedding of one of the two isomorphic components forming the graph $SCay(F_2(\mathrm{var}\{m2_2\}),\{a,b\})$ [right].}
\end{center}

$m3$:

The planar embedding of the graph $SCay(F_2(\mathrm{var}\{m3_1\}),\{a,b\})$ is shown in Fig. 3.

\begin{center}
\includegraphics{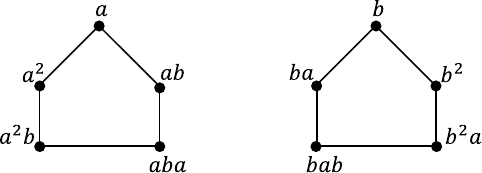}\\
{\footnotesize Fig. 3. Planar embedding of the graph $SCay(F_2(\mathrm{var}\{m3_1\}),\{a,b\})$.}
\end{center}

The graph $SCay(F_3(\mathrm{var}\{m3_1\}),\{a,b,c\})$ contains a subgraph homeomorphic to the graph $K_{3,3}$, since the following pairwise disjoint routes are present between the vertices of the sets $\{a^2,a^2bc,aba\}$ and $\{ab,a^2b,a^2cb\}$: $a^2-a-ab$; $a^2-a^2b$; $a^2-a^2c-a^2cb$; $a^2bc-abca-abc-ab$; $a^2bc-a^2b$; $a^2bc-a^2cb$; $aba-ab$; $aba-a^2b$; $aba-abac-a^2cb$.

The graph $SCay(F_2(\mathrm{var}\{m3_2\}),\{a,b\})$ contains paths between vertices from $\{a^2,aba,ab^2a\}$ and $\{a^3b,a^2b,ab\}$: $a^2-a^3-a^3b$; $a^2-a^2b$; $a^2-a-ab$; $aba-a^3b$; $aba-{\left(ab\right)}^2-a^2ba^2-a^2ba-a^2b$; $aba-ab$; $ab^2a-a^3b^2-a^3b$; $ab^2a-a^2b$; $ab^2a-ab^2-ab$. In the graph $SCay(F_2(\mathrm{var}\{m3_3\}),\{a,b\})$, there are paths between vertices from the sets $\{a^2,aba,a{b}^{2}\}$ and $\{a^{4}b,a^2b,ab\}$: $a^2-a^{4}-a^{4}b$; $a^2-a^2b$; $a^2-a-ab$; $aba-a^{4}b$; $aba-{\left(ab\right)}^2-{\left(ab\right)}^{2}a-{\left(ab\right)}^{2}{a} ^{2}-{a}^{3}{b}^{2}ab-{a}^{3}{b}^ {2}a-{a}^{3}{b}^{2}-{a}^{3}b-{a}^{2}{b}^{3}a-{a}^{2}{b}^{3}-a^2b$; $aba-ab$; $a{b}^{2}-a{b}^{2}a-a^{4}b^2-a^{4}b$; $a{b}^{2}-a{b}^{3}-a{b}^{3}a-a^2b$; $a{b}^{2}-ab$, with changes highlighted. In the general case, in the graph $SCay(F_2(\mathrm{var}\{m3_n\}),\{a,b\})$, where $n>3$, we find paths between vertices from the sets $\{a^2,aba,ab^2\}$ and $\{a^{n+1}b,a^2b,ab\}$: $a^2-a^{n+1}-a^{n+1}b$; $a^2-a^2b$; $a^2-a-ab$; $aba-a^{n+1}b$; $aba-{\left(ab\right)}^2-{\left(ab\right)}^2a-{\left(ab\right)}^2a^2-\dots -a^2b$; $aba-ab$; $ab^2-\dots -a^{n+1}b^2-a^{n+1}b$; $ab^2-ab^3-\dots -a^2b$; $ab^2-ab$.

\begin{center}
\includegraphics{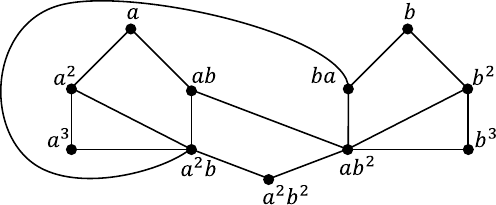}\\
{\footnotesize Fig. 4. Planar embedding of the graph $SCay(F_2(\mathrm{var}\{m4_{1,(123)}\}),\{a,b\})$.}
\end{center}

$m4$:

The planar embedding of the graph $SCay(F_2(\mathrm{var}\{m4_{1,(123)}\}),\{a,b\})$ is shown in Fig. 4.

The existence of a subgraph homeomorphic to the graph $K_{3,3}$ in the simplified Cayley graph of the 3-generated semigroup $F_3(\mathrm{var}\left\{m4_{1,\left(123\right)}\right\})$ is determined by the existence of the following pairwise disjoint routes between the vertices of the set $\{a,b,c\}$ and $\{a^2,b^2,c^2\}$: $a-a^2$; $a-ab-ab^2-b^2$; $a-ac-ac^2-c^2$; $b-ba-a^2b-a^2$; $b-b^2$; $b-bc-bc^2-c^2$; $c-ca-a^2c-a^2$; $c-cb-b^2c-b^2$; $c-c^2.$

At the core of the Cayley graph of the semigroup $F_2(\mathrm{var}\left\{m4_{2,\left(123\right)}\right\})$, a subgraph homeomorphic to the graph $K_{3,3}$ is formed by routes between vertices from the sets $\{a^2,a^2b^2,ba\}$ and $\{a^3b^2,a^2b,ab^2\}$: $a^2-a^3-a^3b-a^3b^2$; $a^2-a^2b$; $a^2-a-ab-ab^2$; $a^2b^2-a^3b^2$; $a^2b^2-a^2b$; $a^2b^2-ab^2$; $ba-b-b^2-b^3-ab^3-a^2b^3-a^3b^3-a^3b^2$; $ba-a^2b$; $ba-ab^2$.

At the core of the Cayley graph of the semigroup $F_2(\mathrm{var}\left\{m4_{3,\left(123\right)}\right\})$, a subgraph homeomorphic to the graph $K_{3,3}$ is formed by routes between vertices from the sets $\{a^2,a^4,a^3b\ \}$ and $\{a^3,a^2b,a^4b\}$: $a^2-a^3$; $a^2-a^2b$; $a^2-a-ab-ab^2-a^2b^2-a^4b^2-a^4b$; $a^4-a^3$; $a^4-a^5-a^2b$; $a^4-a^4b$; $a^3b-a^3$; $a^3b-a^2b$; $a^3b-a^4b$. Similarly, in the simplified Cayley graph of the semigroup $F_2(\mathrm{var}\left\{m4_{4,\left(123\right)}\right\})$, the subgraph homeomorphic to the graph $K_{3,3}$ is formed by routes between vertices from the sets $\{a^2,a^4,a^3b\}$ and $\{a^3,a^2b;a^4b\}$: $a^2-a^3$; $a^2-a^2b$; $a^2-a-ab-ab^2-a^2b^2-{a}^{3}{b}^{2}-a^4b^2-a^4b$; $a^4-a^3$; $a^4-a^5-{a}^{6}-a^2b$; $a^4-a^4b$; $a^3b-a^3$; $a^3b-a^2b$; $a^3b-a^4b$, with different elements highlighted. For each $n>4$ in the simplified Cayley graph of the semigroup $F_2(\mathrm{var}\left\{m4_{n,\left(123\right)}\right\})$ a subgraph homeomorphic to the graph $K_{3,3}$ is formed by routes between vertices from the sets $\{a^2,a^4,a^3b\}$ and $\{a^3,a^2b;a^4b\}$: $a^2-a^3$; $a^2-a^2b$; $a^2-a-ab-ab^2-a^2b^2-a^3b^2-\dots -a^4b^2-a^4b$; $a^4-a^3$; $a^4-a^5-a^6-\dots -a^2b$; $a^4-a^4b$; $a^3b-a^3$; $a^3b-a^2b$; $a^3b-a^4b$.

For any permutation $\pi \in {\prod}_1$, the simplified Cayley graph of the semigroup $F_2(\mathrm{var}\{m4_{1,\pi }\})$ is isomorphic to $SCay(F_2(\mathrm{var}\{m4_{1,(123)}\}),\{a,b\})$, the planar embedding of which is shown in Fig. 4. And for $n>1$, the same forbidden configurations are found in $F_2(\mathrm{var}\{m4_{n,\pi }\})$ as in $F_2(\mathrm{var}\{m4_{n,(123)}\})$.

$m5$:

Moreover, for any permutation $\pi \in {\prod}_1$, the Cayley graph of the semigroup $F_2(\mathrm{var}\left\{m5_{\pi }\right\})$ has a homeomorphic $K_{3,3}$ subgraph on pairwise disjoint paths between vertices from the sets $\{a^2,ba,a^2b^2\}$ and $\{ab^2,a^2b,a^4b\}$: $a^2-a-ab-ab^2$; $a^2-a^2b$; $a^2-a^3-a^4-a^4b$; $ba-ab^2$; $ba-a^2b$; $ba-b-b^2-b^3-ab^3-a^4b$; $a^2b^2-ab^2$; $a^2b^2-a^2b$; $a^2b^2-a^4b$.

$m6$:

Similarly, for any permutation $\pi \in {\prod}_1$ in the simplified Cayley graph of the semigroup $F_2(\mathrm{var}\left\{m6_{\pi }\right\})$ there is a homeomorphic $K_{3,3}$ subgraph on pairwise disjoint routes between vertices from the sets $\{a^2,ba,a^2b^2\}$ and $\{ab^2,a^2b,a^2b^3\}$: $a^2-a-ab-ab^2$; $a^2-a^2b$; $a^2-a^3-a^4-a^5-a^5b$; $ba-ab^2$; $ba-a^2b$; $ba-b-b^2-b^3-ab^3-a^2b^3$; $a^2b^2-ab^2$; $a^2b^2-a^2b$; $a^2b^2-a^2b^3$.

$m7$:

For any permutation $\pi \in {\prod}_1$, the Cayley graph of the semigroup $F_2(\mathrm{var}\left\{m7_{\pi }\right\})$ contains a homeomorphic $K_{3,3}$ subgraph on pairwise disjoint paths between vertices from the sets $\{a^2,ba,a^2b^2\}$ and $\{ab^2,a^2b,a^4b\}$: $a^2-a-ab-ab^2$; $a^2-a^2b$; $a^2-a^3-a^4-a^4b$; $ba-ab^2$; $ba-a^2b$; $ba-b-b^2-b^3-ab^3-a^4b$; $a^2b^2-ab^2$; $a^2b^2-a^2b$; $a^2b^2-a^4b$.

$m8$:

A planar embedding of the simplified Cayley graph of the semigroup $F_2(\mathrm{var}\left\{m8_{\pi }\right\})$, for any $\pi \in {\prod}_1$, is shown in Fig. 5.

\begin{center}
\includegraphics{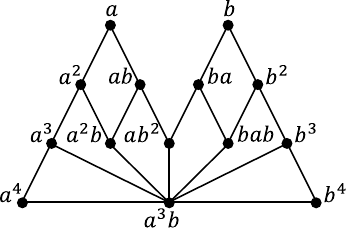}\\
{\footnotesize Fig. 5. Planar embedding of the graph $SCay(F_2(\mathrm{var}\{m8_{\pi }\}),\{a,b\})$, for any $\pi \in {\prod}_1$.}
\end{center}

With further increase in the number of generators in the graph $SCay(F_3(\mathrm{var}\{m8_{\pi }\}),\{a,b,c\})$, for $\pi \in {\prod}_1 \setminus \{(234)\}$, a subgraph homeomorphic to the graph $K_{3,3}$ is found on the following routes between vertices from $\{a^2,a^2b,ab^2\}$ and $\{a^3b,a^2bc,ab\}$: $a^2-a^3-a^3b$; $a^2-a^2c-a^2bc$; $a^2-a-ab$; $a^2b-a^3b$; $a^2b-a^2bc$; $a^2b-ab$; $ab^2-a^3b$; $ab^2-a^2bc$; $ab^2-ab$. For the excluded case, in the graph $SCay(F_3(\mathrm{var}\{m8_{(234)}\}),\{a,b,c\})$, it is easy to find a subgraph homeomorphic to the graph $K_{3,3}$ on the following routes between vertices from $\{a^2,ab,ac\}$ and $\{a,a^2b,ab^2\}$: $a^2-a$; $a^2-a^2b$; $a^2-a^3-a^3b-ab^2$; $ab-a$; $ab-a^2b$; $ab-ab^2$; $ac-a$; $ac-a^2c-a^2bc-a^2b$; $ac-ac^2-ca-c-cb-bc^2-bc-b-ba-ab^2$.

$m9$:

For any permutation $\pi \in {\prod}_1$, the Cayley graph of the semigroup $F_2(\mathrm{var}\left\{m9_{\pi }\right\})$ contains a homeomorphic $K_{3,3}$ subgraph on pairwise disjoint paths between vertices from the sets $\{ab,ba,a^3b\}$ and $\{aba,ab^2,a^2b\}$: $ab-aba$; $ab-ab^2$; $ab-a-a^2-a^2b$; $ba-aba$; $ba-ab^2$; $ba-b-b^2-a^2b$; $a^3b-aba$; $a^3b-ab^2$; $a^3b-a^2b$.

$m10$:

For any permutation $\pi \in {\prod}_1$ in the simplified Cayley graph of the semigroup $F_2(\mathrm{var}\left\{m{10}_{\pi }\right\})$ there is a homeomorphic $K_{3,3}$ subgraph on pairwise disjoint routes between vertices from the sets $\{ab,ba,a^3b\}$ and $\{aba,ab^2,a^2b\} $: $ab-aba$; $ab-ab^2$; $ab-a-a^2-a^2b$; $ba-b-b^2-aba$; $ba-ab^2$; $ba-a ^2b$; $a^3b-aba$; $a^3b-ab^2$; $a^3b-a^2b$.

$m11$:

A planar embedding of the simplified Cayley graph of the semigroup $F_2(\mathrm{var}\left\{m{11}_{\pi }\right\})$, for any $\pi \in {\prod}_1$, is shown in Fig.6.

\begin{center}
\includegraphics{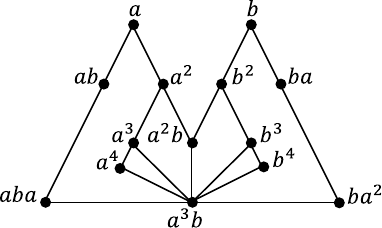}\\
{\footnotesize Fig.6. Planar embedding of the graph $SCay(F_2(\mathrm{var}\{m{11}_{\pi }\}),\{a,b\})$, for any $\pi \in {\prod}_1$.}
\end{center}

And in the column $SCay(F_3(\mathrm{var}\{m{11}_{\pi }\}),\{a,b,c \})$, for any $\pi \in {\prod}_1$, a subgraph homeomorphic to the graph $K_{3,3}$ is found on the following routes between vertices that are elements of the set $\{a,a ^3,a^2b\}$ and $\{a^2,a^3b,a^3c\}$: $a-a^2$; $a-ab-aba-a^3b$; $a-ac-aca-a^3c$; $a^3-a^2$; $a^3-a^3b$; $a^3-a^3c$; $a^2b-a^2 $; $a^2b-a^3b$; $a^2b-b^2-b^2c-c^2-a^2c-a^3c$.

$m12$:

Planar embedding of the simplified Cayley graph of the semigroup $F_2 (\mathrm{var}\left\{m{12}_{\pi }\right\})$, for any $\pi \in {\prod}_1$, is shown in Fig. 7.

\begin{center}
\includegraphics{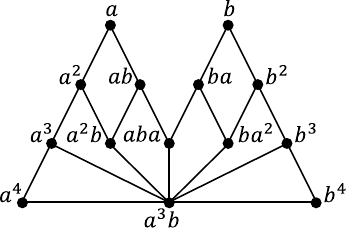}\\
{\footnotesize Fig.7. Planar embedding of the graph $SCay(F_2(\mathrm{var}\{m{12}_{\pi }\}),\{a,b\})$, for any $\pi \in {\prod}_1$.}
\end{center}

In the column $SCay(F_3(\mathrm{var}\{m{12}_{\pi }\}),\{a,b,c\})$, for each $\pi \in \left\{\left(123\right),\ \left(124\right),\ \left(134\right),\ \left(234\right)\right\}\?={\prod}_1\backslash \{\left(12\right)\left(34\right),\left(13\right)\left(24\right),\left(14\right)\left(23\right)\}$, a subgraph homeomorphic to the graph $K_{3,3}$ is found on the following routes between vertices from $\{a^2,aba,a^2b\}$ and $\{ab,a^3b,a^2bc\}$: $a^2-a-ab$; $a^2-a^3-a^3b$; $a^2-a^2c-a^2bc$; $aba-ab$; $aba-a^3b $; $aba-a^2bc$; $a^2b-ab$; $a^2b-a^3b$; $a^2b-a^2bc$. And for the remaining $\pi \in \{\left(12\right)\left(34\right),\left(13\right)\left(24\right),\left(14\right)\left( 23\right)\}$ it is easy to discover routes between vertices from $\{a^2,ab,ac\}$ and $\{a,a^2b,aba\}$: $a^2-a$; $a^2-a^2b$; $a^2-a^3-a^3b-aba$; $ab-a$; $ab-a^2b$; $ab-aba$; $ac-a$; $ac-a^2c-a^2bc-a^2b$; $ac-aca-ca-c-cb-bcb-bc-b-ba-aba$.

$m13$:

A planar embedding of the simplified Cayley graph of the semigroup $F_2(\mathrm{var}\left\{m{13}_{\pi }\right\})$, for any $\pi \in {\prod}_1$, is shown in Fig. 8. And already in the graph $SCay(F_3(\mathrm{var}\{m{13}_{\pi }\}),\{a,b,c\}) $, for any $\pi \in {\prod}_1$, a subgraph homeomorphic to the graph $K_{3,3}$ is found on the following routes between vertices from $\{a,b,c\}$ and $ \{a^2,b^2,c^2\}$: $a-a^2$; $a-ab-aba-b^2$; $a-ac-aca-c^2$; $b-ba-a^2b-a^2$; $b-b^2$; $b-bc-bcb-c^2$; $c-ca-a^2c-a^2$; $c-cb-b^2c-b^2$; $c-c^2$.

\begin{center}
\includegraphics{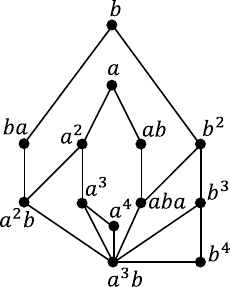}\\
{\footnotesize Fig.8. Planar embedding of the graph $SCay(F_2(\mathrm{var}\{m{13}_{\pi }\}),\{a,b\})$, for any $\pi \in { \textrm{P}}_1$.}
\end{center}

$m14$:

The semigroup $F_3(\mathrm{var}\left\{m{14}_{\pi }\right\})$ has 34 elements. In order not to clutter the proof, we will not give a planar embedding of the simplified Cayley graph of its, but will use the Colin de Verdi\`ere invariant~\cite{4}. The adjacency matrix of the graph $SCay(F_3(\mathrm{var}\{m{14}_{\ pi }\}),\{a,b,c\})$ is graphically depicted in Fig. 9, where in an array of 34x34 dots, a black dot indicates the presence of an edge, and a white one indicates its absence between the vertices of the corresponding row and column. The Colin de Verdi\`ere invariant is defined in a manner similar to the adjacency matrix. Namely, for a simple graph $G$ without loops and multiple edges on $n$ vertices, the Colin de Verdi\`ere invariant $\mu (G)$ is the largest corank of a square matrix $M$ of order $n$ with real coefficients such that its the off-diagonal elements in row $i$ column $j$ are negative numbers when vertex $i$ is adjacent to vertex $j$ in graph $G$, otherwise these elements are zero. In addition, matrix $M$ must have exactly one negative eigenvalue of multiplicity 1 (i.e., be diagonalizable) and there is no nonzero matrix $X$ with real coefficients such that $MX=O$ and $X_{i,j}=0\ $whenever $i= j$ or $M_{i,j}\neq 0$ (that is, $i$ and $j$ are adjacent). Maximization of corank is achieved by changing the diagonal elements of the matrix $M$ in order to minimize its rank. To achieve the restrictions specified in the definition of the Colin de Verdi\`ere invariant, one can take the adjacency matrix of graph $G$ as a starting point, change the signs of all its elements to the opposite, and at the intersection $ i$ rows $i$ columns diagonal elements put $i$, for $i>1$, and -1 for $i=1$, this will provide further diagonalizing the resulting matrix $M$, that is, it will lead to a single negative eigenvalue of multiplicity 1. After that, minimize the rank of the resulting matrix. The graph $G$ turns out to be planar if and only if $\mu \left(G\right)\le 3$. In the case under consideration, we have $\mu \left(SCay\left(F_3\left(\mathrm{var}\left\{m{14}_{\pi }\right\}\right),\left\{a,b,c\right\}\right)\right)\le 3$.

\begin{center}
\includegraphics*[width=1.42in, height=1.42in]{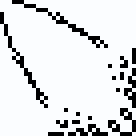}\\
{\footnotesize Fig.9. Adjacency matrix of the graph $SCay(F_3(\mathrm{var}\{m{14}_{\pi }\}),\{a,b,c\})$.}
\end{center}

$m10$:

A planar embedding of the simplified Cayley graph of the semigroup $F_2(\mathrm{var}\left\{m{15}_{\pi }\right\})$, for any $\pi \in {\prod}_1$, is shown in Fig.10. And in the graph $SCay(F_3(\mathrm{var}\{m{15}_{\pi }\}),\{a,b,c\})$, for any $\pi \in {\prod}_1$, the same as in the graph $SCay(F_3(\mathrm{var}\{m{13}_{\pi }\}),\{a,b,c\})$, a subgraph homeomorphic to the graph $K_{3,3}$ on the routes between the vertices from $\{a,b,c\}$ and $\left\{a^2,b^2,c^2\right\}$ is found.

\begin{center}
\includegraphics{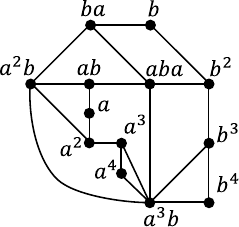}\\
{\footnotesize Fig.10. Planar embedding of the graph $SCay(F_2(\mathrm{var}\{m{15}_{\pi }\}),\{a,b\})$, for any $\pi \in {\prod}_1$.}
\end{center}

$m16$:

A planar embedding of the simplified Cayley graph of the semigroup $F_2(\mathrm{var}\left\{m{16}_{\pi }\right\})$, for any $\pi \in {\prod}_1$, is shown in Fig.11. And in the graph $SCay(F_3(\mathrm{var}\{m{16}_{\pi }\}),\{a,b,c\})$, for any $\pi \in {\prod}_1$, a subgraph homeomorphic to the graph $K_{3,3}$ is found on the following routes between vertices from $\{a,a^3,a^2b\}$ and $\{a^2,aba,aca\}$: $a-a^2$; $a-ab-aba$; $a-ac-aca$; $a^3-a^2$; $a^3-aba$; $a^3-aca$; $a^2b-a^2$; $a^2b-aba$; $a^2b-b^2-b^2c-c^2-c-ca-aca$.

\begin{center}
\includegraphics{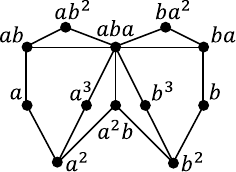}\\
{\footnotesize Fig. 11. Planar embedding of the graph $SCay(F_2(\mathrm{var}\{m{16}_{\pi }\}),\{a,b\})$, for any $\pi \in {\prod}_1$.}
\end{center}

$m17$:

A planar embedding of the simplified Cayley graph of the semigroup $F_2(\mathrm{var}\left\{m{17}_{\pi }\right\})$, for any $\pi \in {\prod}_1$, is shown in Fig. 12.

\begin{center}
\includegraphics{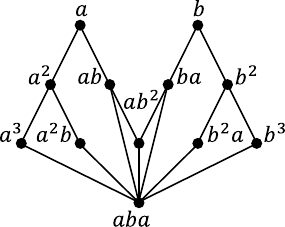}\\
{\footnotesize Fig. 12. Planar embedding of the graph $SCay(F_2(\mathrm{var}\{m{17}_{\pi }\}),\{a,b\})$, for any $\pi \in {\prod}_1$.}
\end{center}

For any $\pi\in\left\{\left(123\right),\ \left(124\right),\ \left(134\right),\ \left(234\right),\left(13\right)\left(24\right)\right\}\?={\prod}_1\backslash \{\left(12\right)\left(34\right),\left(14\right)\left(23\right)\}$ in the graph $SCay(F_3(\mathrm{var}\{m{17}_{\pi }\}),\{a,b,c\})$ found subgraph homeomorphic to the graph $K_{3,3}$ on the following routes between vertices from $\{a^2,aba,aca\}$ and $\{a,a^2b,a^3\}$: $a^2-a$; $a^2-a^2b$; $a^2-a^3$; $aba-ab-a$; $aba-a^2b$; $aba-a^3$; $aca-ac-a$; $aca-a^2bc-a^2b$; $aca-a^3$.

In the graph $SCay(F_3(\mathrm{var}\{m{17}_{\left(12\right)\left(34\right)}\}),\{a,b,c\})$, the presence of a subgraph homeomorphic to the graph $K_{3,3}$ is due to the existence of routes between vertices from $\{a^2,aba,aca\}$ and $\{a,a^2b,a^3\}$: $a^2-a$; $a^2-a^2b$; $a^2-a^3$; $aba-ab-a$; $aba-a^2b$; $aba-a^3$; $aca-ac-a$; $aca-a^2c-a^2bc-a^2b$; $aca-a^3$. And in the remaining graph $SCay(F_3(\mathrm{var}\{m{17}_{(14)(23)}\}),\{a,b,c\})$ the presence of a subgraph homeomorphic to the graph $K_{3,3}$ is due to the existence of routes between the vertices from $\{a^2,aba,aca\}$ and $\{a,a^2c,a^3\}$: $a^2-a$; $a^2-a^2c$; $a^2-a^3$; $aba-ab-a$; $aba-a^2bc-a^2c$; $aba-a^3$; $aca-ac-a$; $aca-a^2c$; $aca-a^3$.

$m18$:

A planar embedding of the simplified Cayley graph of the semigroup $F_2(\mathrm{var}\left\{m{18}_{\pi }\right\})$, for any $\pi \in {\prod}_1$, is shown in Fig. 13.

\begin{center}
\includegraphics{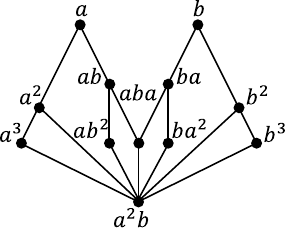}\\
{\footnotesize Fig. 13. Planar embedding of the graph $SCay(F_2(\mathrm{var}\{m{18}_{\pi }\}),\{a,b\})$, for any $\pi \in {\prod}_1$.}
\end{center}

In the graph $SCay(F_3(\mathrm{var}\{m{18}_{\pi }\}),\{a,b,c\})$, for each $\pi \in \left\{\left(123\right),\ \left(124\right),\ \left(134\right),\ \left(234\right)\right\}\?={\prod}_1\backslash \{\left(12\right)\left(34\right),\left(13\right)\left(24\right),\left(14\right)\left(23\right)\}$, a subgraph homeomorphic to the graph $K_{3,3}$ is found on the following routes between vertices from $\{a^2,aba,ab^2\}$ and $\{ab,a^2b,a^2bc\}$: $a^2-a-ab$; $a^2-a^2b$; $a^2-a^2c-a^2bc$; $aba-ab$; $aba-a^2b$; $aba-a^2bc$; $ab^2-ab$; $ab^2-a^2b$; $ab^2-a^2bc$. And with the remaining $\pi \in \{\left(12\right)\left(34\right),\left(13\right)\left(24\right),\left(14\right)\left(23\right)\}$, we can find routes between vertices from $\{a^2,b^2,a^2bc\}$ and $\{ab,a^2b,a^2c\}$: $a^2-a-ab$; $a^2-a^2b$; $a^2-a^2c$; $b^2-b-ba-aba-ab$; $b^2-a^2b$; $b^2-b^2c-c^2-a^2c$; $a^2bc-ab^2-ab$; $a^2bc-a^2b$; $a^2bc-a^2c$.

$m19$:

A planar embedding of the simplified Cayley graph of the semigroup $F_2(\mathrm{var}\left\{m{19}_{\pi }\right\})$, for any $\pi \in {\prod}_1$, is shown in Fig. 14.

\begin{center}
\includegraphics{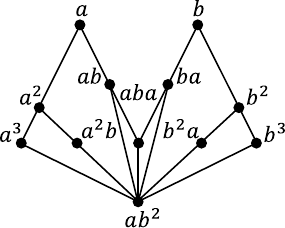}\\
{\footnotesize Fig. 14. Planar embedding of the graph $SCay(F_2(\mathrm{var}\{m{19}_{\pi }\}),\{a,b\})$, for any $\pi \in {\prod}_1$.}
\end{center}

In the graph $SCay(F_3(\mathrm{var}\{m{19}_{\pi }\}),\{a,b,c\})$, for each $\pi \in \left\{\left(123\right),\ \left(124\right),\ \left(134\right),\ \left(234\right)\right\}\?={\prod}_1\backslash \{\left(12\right)\left(34\right),\left(13\right)\left(24\right),\left(14\right)\left(23\right)\}$, a subgraph homeomorphic to the graph $K_{3,3}$ is found on the following routes between vertices from $\{a^2,ab^2,ac^2\}$ and $\{a,a^3,a^2b\}$: $a^2-a$; $a^2-a^3$; $a^2-a^2b$; $ab^2-ab-a$; $ab^2-a^3$; $ab^2-a^2b$; $ac^2-ac-a$; $ac^2-a^3$; $ac^2-a^2bc-a^2b$. For $\pi \in \{\left(12\right)\left(34\right),\left(13\right)\left(24\right)\}$, on the following routes between vertices from $\{a^2,ab^2,ac^2\}$ and $\{a,a^3,a^2b\}$: $a^2-a$; $a^2-a^3$; $a^2-a^2b$; $ab^2-ab-a$; $ab^2-a^3$; $ab^2-a^2b$; $ac^2-ac-a$; $ac^2-a^3$; $ac^2-a^2c-a^2bc-a^2b$. And with the remaining $\pi =\left(14\right)\left(23\right)$, on the routes between the vertices from $\{a^2,ab^2,ac^2\}$ and $\{a,a^3,a^2c\}$: $a^2-a$; $a^2-a^3$; $a^2-a^2c$; $ab^2-ab-a$; $ab^2-a^3$; $ab^2-a^2bc-a^2c$; $ac^2-ac-a$; $ac^2-a^3$; $ac^2-a^2c$.

$m20\div 23$:

In the graphs $SCay\left(F_2\left(\mathrm{var}\left\{m{20}_{\pi }\right\}\right),\left\{a,b\right\}\right)$, for any $\pi \in {\textrm{P}}_2$, $SCay\left(F_2\left(\mathrm{var}\left\{m21\right\}\right),\left\{a,b\right\}\right)$, $SCay\left(F_2\left(\mathrm{var}\left\{m22\right\}\right),\left\{a,b\right\}\right)$ and $SCay\left(F_2\left(\mathrm{var}\left\{m23\right\}\right),\left\{a,b\right\}\right)$, a subgraph homeomorphic to the graph $K_{3,3}$ is found on the following routes between vertices from the set $\left\{a^2b,aba,ab^2\right\}$ and $\left\{ba,ab,a^3b\right\}$: $a^2b-ba$; $a^2b-a^2-a-ab$; $a^2b-a^3b$; $aba-ba$; $aba-ab$; $aba-a^3b$; $ab^2-b^2-b-ba$; $ab^2-ab$; $ab^2-a^3b$.

$m24$:

A planar embedding of the simplified Cayley graph of the semigroup $F_2(\mathrm{var}\left\{m{24}_{\pi }\right\})$, for any permutation $\pi \in {\prod}_3\backslash \{\left(12\right)\left(34\right)\}$, is shown in Fig. 15 on the left. Note that it is isomorphic to the simplified Cayley graph of the semigroup $F_2(\mathrm{var}\left\{m{18}_{\pi }\right\})$, for any $\pi \in {\prod}_1$, shown in Fig. 13. And the planar embedding of the simplified Cayley graph of the semigroup $F_2(\mathrm{var}\left\{m{24}_{(12)(34)}\right\})$ is shown in Fig.15 on the right.

\begin{center}
\includegraphics{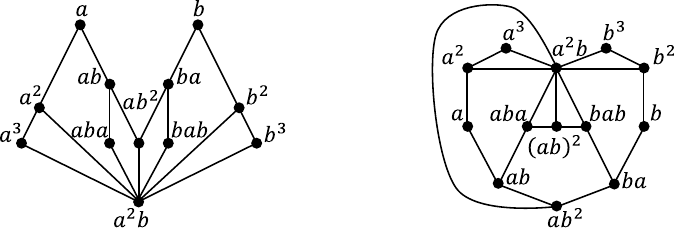}\\
{\footnotesize Fig.15. Planar embedding of the graph $SCay\left(F_2\left(\mathrm{var}\left\{{\sigma }_1\right\}\right),\left\{a,b\right\}\right)$, \\where ${\sigma }_1\in \{m{24}_{\pi },m27,m28,m32\}$, for any $\pi \in {\prod}_3\backslash \{\left(12\right)\left(34\right)\}$ [on the left] and\\ $SCay\left(F_2\left(\mathrm{var}\left\{{\sigma }_2\right\}\right),\left\{a,b\right\}\right)$, where ${\sigma }_2\in \{m{24}_{\left(12\right)\left(34\right)},m26,m36\}$ [right].}
\end{center}

In the $SCay\left(F_3\left(\mathrm{var}\left\{m{24}_{\pi }\right\}\right),\left\{a,b,c\right\}\right)$, for $\pi \in {\prod}_3\backslash \{\left(12\right)\left(34\right)\}$, by increasing the number of generators we find a subgraph homeomorphic to the graph $K_{3,3}$ on the following routes between vertices from $\left\{a^2,aba,ab^2\right\}$ and $\left\{ab,a^2b,a^2bc\right\}$: $a^2-a-ab$; $a^2-a^2b$; $a^2-a^2c-a^2bc$; $aba-ab$; $aba-a^2b$; $aba-a^2bc$; $ab^2-ab$; $ab^2-a^2b$; $ab^2-a^2bc$. And in the remaining graph $SCay\left(F_3\left(\mathrm{var}\left\{m{24}_{\{(12)(34)\}}\right\}\right),\left\{a,b,c\right\}\right)$, a subgraph homeomorphic to the graph $K_{3,3}$ exists on the following routes between vertices from $\left\{a^2,b^2,a^2bc\right\}$ and $\left\{ab^2,a^2b,a^2c\right\}$: $a^2-a-ab-ab^2$; $a^2-a^2b$; $a^2-a^2c$; $b^2-b-ba-ab^2$; $b^2-a^2b$; $b^2-b^2c-c^2-a^2c$; $a^2bc-ab^2$; $a^2bc-a^2b$; $a^2bc-a^2c$.

$m25$:

A planar embedding of the simplified Cayley graph of the semigroup $F_2(\mathrm{var}\left\{m{25}_{\pi }\right\})$, with $\pi \in {\prod}_3\backslash \{\left(12\right)\left(34\right)\}$, is shown in Fig.16 on the left. Note that, as in the case of $F_2(\mathrm{var}\left\{m{24}_{\pi }\right\})$, it is isomorphic to the simplified Cayley graph of the semigroup $F_2(\mathrm{var}\left\{m{18}_{\pi }\right\})$, for any $\pi \?\in {\prod}_1$, shown in Fig. 13. And the planar embedding of the simplified Cayley graph of the semigroup $F_2(\mathrm{var}\left\{m{25}_{(12)(34)}\right\})$ is shown in Fig. 16 on the right.

\begin{center}
\includegraphics{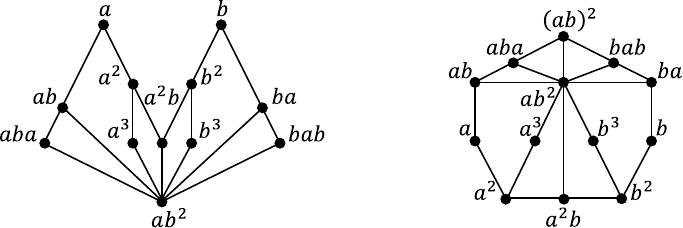}\\
{\footnotesize Fig. 16. Planar embedding of the graph $SCay(F_2(\mathrm{var}\{{\sigma }_1\}),\{a,b\})$, for ${\sigma }_1\in \{m{25}_{\pi },m30,m31,m37,m38,m39,m40\}$, for any $\pi \in {\prod}_3\backslash \{\left(12\right)\left(34\right)\}$ [left] and $SCay\left(F_2\left(\mathrm{var}\left\{{\sigma }_2\right\}\right),\left\{a,b\right\}\right)$, for ${\sigma }_2\in \{m{25}_{\left(12\right)\left(34\right)},m29,m41\}$ [right].}
\end{center}

By increasing the number of generators, in the graphs $SCay\left(F_3\left(\mathrm{var}\left\{m{25}_{\pi }\right\}\right),\left\{a,b,c\right\}\right)$, for any $\pi \in {\prod}_3\backslash \{\left(12\right)\left(34\right)\}$, and in the graph $SCay\left(F_3\left(\mathrm{var}\left\{m{25}_{\{(12)(34)\}}\right\}\right),\left\{a,b,c\right\}\right)$, which is not isomorphic to them, we find the same subgraph homeomorphic to the graph $K_{3,3}$ on the following routes between vertices from $\left\{a^2,ab^2,ac^2\right\}$ and $\left\{a,a^2b,a^3\right\}$: $a^2-a$; $a^2-a^2b$; $a^2-a^3$; $ab^2-ab-a$; $ab^2-a^2b$; $ab^2-a^3$; $ac^2-ac-a$; $ac^2-ca-c-c^2-b^2c-b^2-a^2b$; $ac^2-a^3$.

$m26$:

The planar embedding of the simplified Cayley graph of the semigroup $F_2(\mathrm{var}\left\{m26\right\})$, which coincides with the planar embedding of the simplified Cayley graph of the semigroup $F_2(\mathrm{var}\left\{m{24}_{(12)(34)}\right\})$, is shown in Fig.15 on the right. And in $SCay\left(F_3\left(\mathrm{var}\left\{m26\right\}\right),\left\{a,b,c\right\}\right)$ a subgraph homeomorphic to the graph $K_{3,3}$ is found, the same as in the simplified Cayley graph $SCay\left(F_3\left(\mathrm{var}\left\{m{24}_{\{(12)(34)\}}\right\}\right),\left\{a,b,c\right\}\right)$, on the routes between the vertices from $\left\{a^2,b^2,a^2bc\right\}$ and $\left\{ab^2,a^2b,a^2c\right\}$.

$m27\div 28$:

The planar embedding of the simplified Cayley graph of the semigroups $F_2(\mathrm{var}\left\{m27\right\})$ and $F_2(\mathrm{var}\left\{m28\right\})$, which coincides with the planar embedding of the simplified Cayley graph of the semigroup $F_2(\mathrm{var}\left\{m{24}_{\pi }\right\})$, for any permutation of $\pi \in {\prod}_3\backslash \{\left(12\right)\left(34\right)\}$, is shown in Fig. 15 on the left. And in $SCay\left(F_3\left(\mathrm{var}\left\{m27\right\}\right),\left\{a,b,c\right\}\right)$ with $SCay\left(F_3\left(\mathrm{var}\left\{m28\right\}\right),\left\{a,b,c\right\}\right)$ a subgraph homeomorphic to the graph $K_{3,3}$ is found, the same as in the simplified Cayley graph $SCay\left(F_3\left(\mathrm{var}\left\{m{24}_{\{(12)(34)\}}\right\}\right),\left\{a,b,c\right\}\right)$, on the routes between the vertices from $\left\{a^2,b^2,a^2bc\right\}$ and $\left\{ab^2,a^2b,a^2c\right\}$.

$m29$:

A planar embedding of the simplified Cayley graph of the semigroup $F_2(\mathrm{var}\left\{m29\right\})$, which coincides with the planar embedding of the simplified Cayley graph of the $F_2\left(\mathrm{var}\left\{m{25}_{\left(12\right)\left(34\right)}\right\}\right)$ is shown in Fig.16 on the right. And in $SCay\left(F_3\left(\mathrm{var}\left\{m29\right\}\right),\left\{a,b,c\right\}\right)$ a subgraph homeomorphic to the graph $K_{3,3}$ is found, the same as in the $SCay\left(F_3\left(\mathrm{var}\left\{m{25}_{\{(12)(34)\}}\right\}\right),\left\{a,b,c\right\}\right)$, on the routes between the vertices from $\left\{a^2,ab^2,ac^2\right\}$ and $\left\{a,a^2b,a^3\right\}$.

$m30\div 31$:

The planar embedding of the simplified Cayley graph of the semigroups $F_2(\mathrm{var}\left\{m30\right\})$ and $F_2(\mathrm{var}\left\{m31\right\})$, which coincides with the planar embedding of the simplified Cayley graph of the semigroup $F_2\left(\mathrm{var}\left\{m{25}_{\pi }\right\}\right)$, for each permutation $\pi \in {\prod}_3\backslash \{\left(12\right)\left(34\right)\}$, is shown in Fig. 16 on the left. And in $SCay\left(F_3\left(\mathrm{var}\left\{m30\right\}\right),\left\{a,b,c\right\}\right)$ and $SCay\left(F_3\left(\mathrm{var}\left\{m31\right\}\right),\left\{a,b,c\right\}\right)$ a subgraph homeomorphic to the graph $K_{3,3}$ is found, the same as in the simplified Cayley graph $SCay\left(F_3\left(\mathrm{var}\left\{m{25}_{\{(12)(34)\}}\right\}\right),\left\{a,b,c\right\}\right)$, on the routes between the vertices from $\left\{a^2,ab^2,ac^2\right\}$ and $\left\{a,a^2b,a^3\right\}$.

$m32\div 35$:

For all $\sigma \in \{m32,m33,m34,m35\}$ the planar embedding of the simplified Cayley graph of the semigroup $F_2(\mathrm{var}\left\{\sigma \right\})$, as well as $F_2(\mathrm{var}\left\{m27\right\})$, $F_2(\mathrm{var}\left\{m28\right\})$, coincides with the planar embedding of the simplified Cayley graph of the semigroup $F_2(\mathrm{var}\left\{m{24}_{\pi }\right\})$, for any $\pi \in {\prod}_3\backslash \{\left(12\right)\left(34\right)\}$, and is shown in Fig.15 on the left. And in $SCay\left(F_3\left(\mathrm{var}\left\{\sigma \right\}\right),\left\{a,b,c\right\}\right)$ a subgraph homeomorphic to the graph $K_{3,3}$ is found, the same as in the simplified Cayley graph $SCay\left(F_3\left(\mathrm{var}\left\{m{24}_{\{(12)(34)\}}\right\}\right),\left\{a,b,c\right\}\right)$, on the routes between the vertices from $\left\{a^2,b^2,a^2bc\right\}$ and $\left\{ab^2,a^2b,a^2c\right\}$.

$m36$:

A planar embedding of the simplified Cayley graph of the semigroup $F_2(\mathrm{var}\left\{m36\right\})$, which coincides with the planar embedding of the simplified Cayley graph of the $F_2(\mathrm{var}\left\{m{24}_{(12)(34)}\right\})$, is shown in Fig. 15 on the right. And in $SCay\left(F_3\left(\mathrm{var}\left\{m36\right\}\right),\left\{a,b,c\right\}\right)$ a subgraph homeomorphic to the graph $K_{3,3}$ is found, the same as in the $SCay\left(F_3\left(\mathrm{var}\left\{m{24}_{\{(12)(34)\}}\right\}\right),\left\{a,b,c\right\}\right)$, on the routes between the vertices from $\left\{a^2,b^2,a^2bc\right\}$ and $\left\{ab^2,a^2b,a^2c\right\}$.

$m37\div 40$:

Planar embedding of the simplified Cayley graph of each of the semigroups $F_2(\mathrm{var}\left\{m37\right\})$, $F_2(\mathrm{var}\left\{m38\right\})$, $F_2(\mathrm{var}\left\{m39\right\})$ and $F_2(\mathrm{var}\left\{m40\right\})$, coinciding with the planar embedding of the simplified Cayley graph of the semigroup $F_2\left(\mathrm{var}\left\{m{25}_{\pi }\right\}\right)$, for $\pi \in {\prod}_3\backslash \{\left(12\right)\left(34\right)\}$, is shown in Fig. 16 on the left. And in the simplified Cayley graph $SCay\left(F_3\left(\mathrm{var}\left\{m37\right\}\right),\left\{a,b,c\right\}\right)$, in the simplified Cayley graph $SCay\left(F_3\left(\mathrm{var}\left\{m38\right\}\right),\left\{a,b,c\right\}\right)$, in the simplified Cayley graph $SCay\left(F_3\left(\mathrm{var}\left\{m39\right\}\right),\left\{a,b,c\right\}\right)$, $SCay\left(F_3\left(\mathrm{var}\left\{m40\right\}\right),\left\{a,b,c\right\}\right)$ is found to be homeomorphic to the graph $K_{3,3}$ the subgraph is the same as in the $SCay\left(F_3\left(\mathrm{var}\left\{m{25}_{\{(12)(34)\}}\right\}\right),\left\{a,b,c\right\}\right)$, on routes between vertices from $\left\{a^2,ab^2,ac^2\right\}$ and $\left\{a,a^2b,a^3\right\}$.

$m41$:

A planar embedding of the simplified Cayley graph of the semigroup $F_2(\mathrm{var}\left\{m41\right\})$, which coincides with the planar embedding of the simplified Cayley graph of the $F_2\left(\mathrm{var}\left\{m{25}_{\left(12\right)\left(34\right)}\right\}\right)$ is shown in Fig. 16 on the right. And in $SCay\left(F_3\left(\mathrm{var}\left\{m41\right\}\right),\left\{a,b,c\right\}\right)$ a subgraph homeomorphic to the graph $K_{3,3}$ is found, the same as in the $SCay\left(F_3\left(\mathrm{var}\left\{m{25}_{\{(12)(34)\}}\right\}\right),\left\{a,b,c\right\}\right)$, on the routes between the vertices from $\left\{a^2,ab^2,ac^2\right\}$ and $\left\{a,a^2b,a^3\right\}$.

$m42, m44, m46$:

The planar embedding of the simplified Cayley graph of the semigroups $F_2(\mathrm{var}\left\{m42\right\})$, $F_2(\mathrm{var}\left\{m44\right\})$ and $F_2(\mathrm{var}\left\{m46\right\})$ is shown in Fig.17. And with an increased number of generators, starting with three, in $SCay\left(F_3\left(\mathrm{var}\left\{m42\right\}\right),\left\{a,b,c\right\}\right)$, $SCay\left(F_3\left(\mathrm{var}\left\{m44\right\}\right),\left\{a,b,c\right\}\right)$ and $SCay\left(F_3\left(\mathrm{var}\left\{m46\right\}\right),\left\{a,b,c\right\}\right)$ a subgraph homeomorphic to the graph $K_{3,3}$ is already found on the following pairwise disjoint routes between the vertices from the sets $\left\{a^2,b^2,a^3b\right\}$ and $\left\{ab^2,a^2b,a^3\right\}$: $a^2-a-ab-ab^2$; $a^2-a^2b$; $a^2-a^3$; $b^2-b-ba-ab^2$; $b^2-a^2b$; $b^2-b^2c-c^2-a^2c-a^3c-a^3$; $a^3b-ab^2$; $a^3b-a^2b$; $a^3b-a^3$.

\begin{center}
\includegraphics{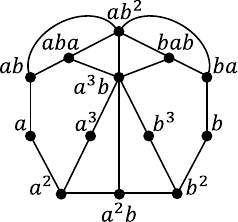}\\
{\footnotesize Fig. 17. Planar embedding of the graph $SCay(F_2(\mathrm{var}\{\sigma \}),\{a,b\})$, for any $\sigma \in \{m42,m44,m46\}$.}
\end{center}

$m43, m45, m47$:

At the simplified Cayley graph of the $F_2\left(\mathrm{var}\left\{m43\right\}\right)$, $F_2\left(\mathrm{var}\left\{m45\right\}\right)$ and $F_2\left(\mathrm{var}\left\{m47\right\}\right)$, a subgraph homeomorphic to the graph $K_{3,3}$ is found on the following pairwise disjoint routes between vertices from the sets $\left\{a^2,aba,ab^2\right\}$ and $\left\{ab,a^2b,a^3b\right\}$: $a^2-a-ab$; $a^2-a^2b$; $a^2-a^3-a^3b$; $aba-ab$; $aba-a^2b$; $aba-a^3b$; $ab^2-ab$; $ab^2-ba-b-b^2-a^2b$; $ab^2-a^3b$.

Thus, the justification of the found values of planarity ranks of the modular varieties of semigroups presented in the theorem being proved was carried out. Theorem 1 is proved.

$\,$
$\,$
\section*{Conclusion}%2

In this paper, all possible modular varieties of semigroups are considered, and the planarity ranks of each of them are calculated. As it turns out, the values of the planarity ranks of an infinitely countable number of all possible modular varieties of semigroups do not exceed 3.

$\,$

$\,$

\end{document}